\documentclass[11pt,leqno,oneside]{amsart}
\usepackage{layout}
\usepackage{mathrsfs,dsfont}
\usepackage{amsmath,amstext,amsthm,amssymb,bbm,color}
\usepackage{comment}
\usepackage{charter}
\usepackage{typearea}
\usepackage{pdfsync}
\usepackage[width=6.5in,height=8.7in]{geometry}

\usepackage{amsfonts}

\DeclareMathOperator\supp{supp}

\def\bt{\begin{thm}}
\def\et{\end{thm}}
\def\bl{\begin{lem}}
\def\el{\end{lem}}
\def\bd{\begin{defn}}
\def\ed{\end{defn}}
\def\bc{\begin{cor}}
\def\ec{\end{cor}}
\def\bp{\begin{proof}}
\def\ep{\end{proof}}
\def\br{\begin{rem}}
\def\er{\end{rem}}

\newtheorem{theorem}{Theorem}[section]
\newtheorem{lemma}[theorem]{Lemma}
\newtheorem{proposition}[theorem]{Proposition}
\theoremstyle{definition}
\newtheorem{definition}[theorem]{Definition}
\newtheorem{example}[theorem]{Example}
\newtheorem{corollary}[theorem]{Corollary}

\theoremstyle{remark}
\newtheorem{remark}[theorem]{Remark}

\numberwithin{equation}{section}

\newcommand{\bthm}{\begin{thm}}
\newcommand{\ethm}{\end{thm}}
\newcommand{\bstp}{\begin{stp}}
\newcommand{\estp}{\end{stp}}
\newcommand{\blemma}{\begin{lemma}}
\newcommand{\elemma}{\end{lemma}}
\newcommand{\bprop}{\begin{prop}}
\newcommand{\eprop}{\end{prop}}
\newcommand{\bpf}{\begin{pf}}
\newcommand{\epf}{\end{pf}}
\newcommand{\bdefn}{\begin{defn}}
\newcommand{\edefn}{\end{defn}}
\newcommand{\brk}{\begin{rmrk}}
\newcommand{\erk}{\end{rmrk}}
\newcommand{\bcrl}{\begin{crl}}
\newcommand{\ecrl}{\end{crl}}

\usepackage{amsfonts}

\newcommand{\overbar}[1]{\mkern 1.5mu\overline{\mkern-1.5mu#1\mkern-1.5mu}\mkern 1.5mu}

\title[$m$-core]{Core sets in K\"{a}hler manifolds  }

\author{N\.Ihat G\"okhan G\"o\u{g}\"u\c{s}, Ozan G\"uny\"uz,  \"{O}zcan Yaz{\i}c{\i}}

\address{Faculty of Engineering and Natural Sciences, Sabancı University, İstanbul, Turkey}

\email{gokhan.gogus@sabanciuniv.edu}
\email{ozangunyuz@alumni.sabanciuniv.edu}

\address{Department  of Mathematics,  Middle East  Technical University,  06800 Ankara,  Turkey }
\email{oyazici@metu.edu.tr}

\date{\today}


\begin{document}

\begin{abstract}
The primary objective of this paper is to study core sets in the setting of $m$-subharmonic functions on the class of (non-compact) K\"{a}hler manifolds. Core sets are investigated in different aspects by considering various classes of plurisubharmonic functions. One of the crucial concepts in studying the structure of this kind of sets is the pseudoconcavity. In a more general way, we will have the structure of core defined with respect to the $m$-subharmonic functions, which we call $m$-core in our setting, in terms of $m$-pseudoconcave sets. In the context of $m$-subharmonic functions, we define $m$-harmonic functions and show that, in $\mathbb{C}^n\,\,(n\geq 2)$ and more generally in any K\"{a}hler manifold of dimension at least $2$, $m$-harmonic functions are pluriharmonic functions for $m \geq 2$.

  \end{abstract}
\maketitle

\section{INTRODUCTION, PRELIMINARIES AND THE RESULT}

The notion core of a complex manifold $M$, denoted by $\mathbf{c}(M)$, was first described and studied by Harz, Tomassini and Shcherbina in a series of papers \cite{HST} and \cite{HST2} for strictly pseudoconvex domains in $\mathbb{C}^{n}$ and then in complex manifolds generally. It is the biggest set on which every bounded and continuous plurisubharmonic function on $M$ fails to be strictly plurisubharmonic. In \cite{PS}, Poletsky and Shcherbina employed a modified definition of core (see section 2 below) and established the structure of core defined this way by decomposing it into the foliation sets on which every plurisubharmonic and upper bounded function becomes constant and answered affirmatively a question posed in \cite{HST}. The same  problem was also solved in \cite{Slod2} with completely different methods. Slodkowski, in his paper \cite{Slod2}, also generalized the core by appealing to the sheaves, for background see section 4 in the said paper.

For various regularity classes, one can define the corresponding cores accordingly. Sets of this sort give a chain of inclusions when defined appropriately according to their regularity properties, however not so much information is available between the relations of the various core sets forming this chain. In connection with this chain of inclusions, in a recent study of Harz (\cite{T}), it has been proven that in the aforementioned chain, the first three inclusions are proper.

We find it worthwhile to mention that more than a decade before the recent studies on the subject core,  the authors Slodkowski and Tomassini introduced in their paper \cite{SlT} a very similar concept, called \textit{minimal kernel}, to the core in the so-called weakly complete complex spaces (i.e. complex spaces having a $\mathcal{C}^{k}$-smooth,\, $k \geq 0$,\, plurisubharmonic exhaustion function), which are the most general framework possible thus far. In the same paper, it turns out that minimal kernels determine how far a given complex space is from being Stein.

In this paper, we generalize the core sets by using the $m$-subharmonic functions. The roots of $m$-subharmonic functions go back to the paper of Li (\cite{Li}) in which he gives the definition of  $m$-pseudoconvexity and has the generalization of the existence of a unique classical solution for the Dirichlet problem of symmetric function of the eigenvalues of real hessian matrix of a function defined on a domain in $\mathbb{R}^{n}$ (which has been proved by Caffarelli, Kohn, Nirenberg and Spruck in \cite{CKNS}) for smoothly bounded domains in $\mathbb{C}^{n}$. The important instrument to study $m$-subharmonic functions in $\mathbb{C}^{n}$ is complex $k$-Hessian equation which, to our knowledge, was also first investigated by Li in his paper \cite{Li}. To mention some references, in the works of Blocki, Kolodziej, Dinew, Sadullaev, Hou, Li, these functions are dealt with in different directions related to the various problems, so for a thorough investigation of $m$-subharmonic functions, the reader can consult the papers such as \cite{Hou}, \cite{Bl} \cite{DK}, \cite{Pl}, \cite{SA} and references therein.

Throughout we will let $M$ be a non-compact K\"{a}hler manifold of dimension $n$ with a fixed K\"{a}hler form $\omega$ on it. All complex manifolds considered here are assumed to be countable at infinity.

Let $\Omega \subset M$ be a domain. $\mathcal{C}^{\infty}_{0}(\Omega)$ will denote the space of test functions on $\Omega$, i.e., of infinitely differentiable functions with some compact support in $\Omega$. Similarly, $\mathcal{D}^{p, q}(M)$ denotes the space of test forms of bidegree $(p, q)$ on the complex manifold $M$ and  we will let $\mathcal{D'}_{p, q}(M)$ denote the space of currents of bidegree $(p, q)$ on $M$, so $\langle T, \varphi\rangle=T(\varphi)$  means the pairing of \,$T\in \mathcal{D'}_{p, q}(M)$ and $\varphi\in \mathcal{D}^{n-p, n-q}(M).$

Our main objects are $m$-subharmonic functions and we will utilize the definitions used in \cite{Pl}, except that the definition of a strictly $m$-subharmonic function differs somewhat, see \cite{HST}. A function $u\in \mathcal{C}^2(\Omega)$ is said to be \textsl{$m$-subharmonic} on $\Omega$ if \begin{equation} \label{msub} (dd^{c}u)^{k} \wedge \omega^{n-k}\geq 0\end{equation} for $k=1, 2, \ldots, m$.

A locally integrable function $u: \Omega \rightarrow [-\infty, \infty)$ is called \textsl{$m$-subharmonic }on $\Omega$ if $u$ is upper semicontinuous and \begin{equation}\label{10} dd^{c}u \wedge dd^{c}u_{1} \wedge \ldots \wedge dd^{c}u_{m-1} \wedge \omega^{n-m} \geq 0\end{equation} holds in the weak sense of currents for any $m$-subharmonic $\mathcal{C}^{2}$ functions $u_{1}, \ldots, u_{m-1}$ defined on $\Omega$. According to these definitions,  $n$-subharmonic functions are plurisubharmonic functions  and $1$-subharmonic functions are subharmonic ones.  We denote the class of all $m$-subharmonic functions on $\Omega$ \,by $\mathrm{SH}_{m}(\Omega)$.

We say that a function $u: \Omega \rightarrow \mathbb{R}$\, is\, \textsl{strictly $m$-subharmonic} on $\Omega$ if for any $\phi \in \mathcal{C}_{0}^{\infty}(\Omega)$, there exists an $\epsilon_{0}>0$ such that for every $\epsilon \in (-\epsilon_{0}, \epsilon_{0})$, $u + \epsilon \phi$ is $m$-subharmonic in the sense of (\ref{10}). For $\mathcal{C}^{2}$ functions, of course, one can use the pointwise condition (\ref{msub}). As is easily seen, since $\epsilon=0$ can be taken, any strictly $m$-subharmonic function is $m$-subharmonic, we do not even need to assume that $u$ is upper semicontinuous. Similar to the $m$-subharmonic functions, if $u\in \mathcal{C}^{2}(\Omega)$, then $u$ is strictly $m$-subharmonic if for any $\phi \in \mathcal{C}_{0}^{\infty}(\Omega)$, there exists an $\epsilon_{0}>0$ such that for every $\epsilon \in (-\epsilon_{0}, \epsilon_{0})$, $$(dd^{c}(u + \epsilon \phi))^{k} \wedge \omega^{n-k}\geq 0.$$

We will mainly concentrate on the class of continuous upper bounded $m$-subharmonic functions on  $M$, denoted by $\mathrm{SH}^{cb}_{m}(M)$.

Because of the local nature of the problems we shall be dealing with in the sequel, we enforce the following condition on the definition of $m$-subharmonicity that assures the local approximation of $u$, which shall be our standing assumption throughout the present work:

\textit{(*)In a neighborhood of every point, there is a decreasing sequence $\{u_{j}\}$ of smooth $m$-subharmonic functions converging to $u$, that is $ (dd^{c}u_{j})^k \wedge \omega^{n-k}\geq 0$.}

 This also indicates why we make the assumption of non-compactness on the K\"{a}hler manifold $M$, because we have merely trivial $m$-subharmonic functions on compact K\"{a}hler manifolds, in which case the approximation condition we have just stated no longer holds.

We give now the definition of the $m$-core of a K\"{a}hler manifold $M$ using $m$-subharmonic functions:
\vspace{5mm}

$\mathbf{c}_{m}(M)=\{z\in M: \text{Every}\,\, \text{function}\,\, \text{of}\,\,  \mathrm{SH}_{m}^{cb}(M)\,\,  \text{fails} \, \text{to} \,\, \text{hold} \,\, \,\text{strict}\,\,  \text{$m$-subharmonicity}\,\,\text{near}\,\,\text{z} \}$.

\vspace{5mm}
The set $\mathbf{c}_{m}(M)$ is closed by its definition. In \cite{PS}, the definition of the core set is a little different, they impose the smoothness condition besides being strictly plurisubharmonic. In the paper \cite{HST}, smooth and plurisubharmonic functions are taken into account: It is the set of all points $z \in M$ such that in a neighborhood of $z$, every smooth $m$-subharmonic function fails to satisfy the strict plurisubharmonicity. In this article, we will be interested in the definition of the core using only the class of continuous (not necessarily smooth, even $\mathcal{C}^{2}$) strictly $m$-subharmonic functions as was done in \cite{Slod2}.

As the main the result of the paper, we prove the following in Section \ref{strum}:

\begin{theorem}\label{decompi}
The core $\mathbf{c}_{m}(M)$ can be disintegrated into $m$-pseudoconcave subsets such that every upper bounded continuous $m$-subharmonic function on $M$ is constant on each of these sets.
\end{theorem}

From Theorem \ref{decompi}, we infer
\begin{corollary}
The $m$-core $\mathbf{c}_{m}(M)$ is empty if and only if the functions in $\mathrm{SH}^{cb}_{m}(M)$ separate the points of $M$.
\end{corollary}
We observe that, as was proved in \cite{PS, Slod2} for the set core in the pluripotential setting, the set $m$-core has the same obstructive nature to separating the points of $M$ by $\mathrm{SH}^{cb}_{m}(M)$.

\section{STRUCTURE OF THE $m$-CORE}\label{strum}
Before examining $m$-pseudoconcavity and $m$-core, we start first with an observation concerned with strictly $m$-subharmonic functions whose easy proof can be done just by the definition of strict $m$-subharmonicity and is left to the reader.

\begin{lemma}\label{strlem}
Let $u_{1}, u_{2}, \ldots, u_{k}\in \mathrm{SH}_{m}(\Omega)$. If, at least, one of $u_{j},\,\,j=1, 2, \ldots, k,$ is strictly $m$-subharmonic, then the sum $\sum_{j=1}^{k}\alpha_{j} u_{j}$ is also strictly $m$-subharmonic for $\alpha_{j}>0,\,j=1, 2, \ldots, k$. \end{lemma}

Let $z_{0}\in M$. A point $z$ is said to belong to $A^{b}_{m}(z_{0})$, respectively to $A^{cb}_{m}(z_{0})$ if $v(z)\leq v(z_{0})$ for any $v\in \mathrm{SH}^{b}_{m}(M)$, respectively for any $v \in \mathrm{SH}^{cb}_{m}(M)$. Obviously we have the inclusion $A^{b}_{m}(z_{0}) \subset A^{cb}_{m}(z_{0})$. We also consider a subset of $A^{cb}_{m}(z_{0})$, denoted by $ (A^{cb}_{m}(z_{0}))_{e}$,  to be defined as the set of elements $z\in M$ so that $v(z)=v(z_{0})$ for all $v\in \mathrm{SH}^{cb}_{m}(M)$. If $A^{cb}_{m}(z_{0})=\{z_{0}\}$, then clearly $A^{cb}_{m}(z_{0})=(A^{cb}_{m}(z_{0}))_{e}$. These sets are simply the $m$-subharmonic analogues of the ones defined and investigated in \cite{PS} and \cite{Slod2}. They are also known as \textit{foliations}. Some of the basic properties that $A^{cb}_{m}(z_{0})$ have are listed below. Proofs of the assertions (1) and (2) are immediate. By using the argument verbatim in the proof of (5) of Proposition 6 and Lemma 4 (which can also be adapted to the $m$-subharmonic setting) in \cite{PS}, one proves (3).

\begin{proposition}
\begin{itemize}
  \item[(1)] $A^{cb}_{m}(z_{0})$ is a closed set.
  \item[(2)] If $z_{1}\in A^{b}_{m}(z_{0})$($z_{1}\in A^{cb}_{m}(z_{0})$), then $A^{b}_{m}(z_{1}) \subset A^{b}_{m}(z_{0})$ ($A^{cb}_{m}(z_{1}) \subset  A^{cb}_{m}(z_{0})$).
  \item[(3)] If $A^{cb}_{m}(z_{0})=\{z_{0}\}$, then $z_{0}\notin \mathbf{c}_{m}(M)$.
\end{itemize}
\end{proposition}

In this paper, we do not focus on the sets $A^{cb}_{m}(z_{0})$ unlike what was done in \cite{PS} because we will follow closely the techniques used in \cite{Slod2}. In \cite{PS}, the authors concentrate first on the sets $A^{cb}(z_{0})$ and obtain the $1$-pseudoconcavity of the sets $A^{cb}(z_{0})$ (Theorem 8 there) based on a strict convexity argument. To do so, given a ($\mathcal{C}^{2}$)strictly plurisubharmonic function with a non-zero differential at a point $0\in  \mathbb{C}^{n}$, they use a well-known theorem that guarantees the existence of a local biholomorphic mapping between the open neighborhoods of $0$ to produce a strictly convex function, for example, see Theorem 2.23 of Chapter 6 in \cite{La}. We do not seem to have such a theorem for $m$-subharmonic functions because, as is well-known, $m$-subharmonic functions are not biholomorphically invariant in general and the Levi form may well be negative. In the context of Theorem 2.23 in \cite{La}, given an $m$-subharmonic function $u$, what one can only say for certain  is that $u \circ h^{-1}$ is subharmonic since the positivity of $dd^c u$ is not impacted by the (local) biholomorphic mapping $h$ used there. Within the proof of the main theorem that we prove in this section, we see that $(A^{cb}_{m}(z_{0}))_{e}$ is $m$-pseudoconcave (see below for the definition). We conjecture that $A^{cb}_{m}(z_{0})$ is also $m$-pseudoconcave.

Let $E$ be a closed set in $M$. We will say \, $E$ \,is \textsl{$m$-pseudoconcave in the sense of Rothstein} if for any $z_{0}\in E$ and for any strictly $m$-subharmonic function $\rho$ defined in a neighborhood $V$ of $z_{0}$, within any relatively compact neighborhood $U\subset V$ with $z_{0}\in U$, there is a point $z\in E\cap U$ where $\rho(z)> \rho(z_{0})$. This is a generalization of $1$-pseudoconcavity in the sense of Rothstein investigated in \cite{PS}. Another important concept related to the $1$-pseudoconcavity in the sense of Rothstein is the local maximum property. For any closed set $C$, being $1$-pseudoconcave in the sense of Rothstein and having the local maximum property are equivalent, see \cite{PS, Slod2}.

As noted from the definition above, they are perfect sets, namely they have no isolated points. We now have another definition which again generalizes the local maximum property.

Let $Y$ be a closed set in $M$. We will say $Y$ has the \textsl{$m$-local maximum property} or is an \textsl{m-local maximum set}  if $Y$ is perfect and for any $w_{0}\in Y$ there is an open neighborhood $U$ of $w_{0}$ in $M$ with compact closure such that if an open set $V\Subset U$ contains $w_{0}$ and the set $K=Y\cap \partial V$ is non-empty, then \begin{equation} \label{mloc}\max_{Y\cap \overline{V}}{u}\leq \max_{K}{u} \end{equation} for any $m$-subharmonic function $u$ on $U$.

As in the plurisubharmonic case, we have the equivalence of $m$-pseudoconcave and $m$-local maximum sets. Proof is identical except that, in the necessity part, one uses Lemma \ref{strict}.
\begin{lemma} \label{maxloc}
Let $M$ be a K\"{a}hler manifold. A closed set $Y \subset M$ has the $m$-local maximum property if and only if it is $m$-pseudoconcave in the sense of Rothstein.
\end{lemma}

We will call $u \in \mathrm{SH}_{m}(M)$ \textsl{m-maximal} if, for any relatively compact domain $D \subset M$ and any $v\in \mathrm{SH}_{m}(D)$ which is upper semicontinuous on $\overbar{D}$, $u\geq v$ on $\partial D$, then $u \geq v$ on $D$.

The following proposition is a simple consequence of definition of an $m$-maximal function which was used in \cite{Slod2} without proof. Since it is essential in the sequel, we give its basic proof for $m$-subharmonic functions.

\begin{proposition}\label{maxx}
Let $\varphi$ be an $m$-subharmonic function on a domain $V \subset M$ and $\psi$ is an $m$-maximal continuous function on $V$. Then for every ball $B\subset V$ with $\overline{B}\subset V$, \begin{equation}\label{msc}\sup_{B}{(\varphi- \psi)} \leq \max_{\partial B}{(\varphi- \psi)}.\end{equation}
\end{proposition}

\begin{proof}
Let a ball $B$ be given with $\overbar{B}\subset V$. Since $\partial B$ is compact, we can write $d= \max_{\partial B}{(\varphi -\psi)}$, so\, $\varphi -\psi \leq d$\, on\, $\partial B$. This last inequality can be written as\, $(\varphi - d)-\psi \leq 0$\, on\, $\partial B$. Now define $\varphi'= \varphi - d$, which is also $m$-subharmonic on V, and in particular, is upper semicontinuous on $\overbar{B}$. Then we have $\varphi' - \psi \leq 0$ on $\partial B$, that is, $\varphi' \leq \psi$ on $\partial B$, but by assumption that $\psi$ is $m$-maximal, one gets that $\varphi' \leq \psi$ on $B$, which is equivalent to $\varphi - \psi \leq d$ on $B$. Taking supremum of the left side, the desired inequality (\ref{msc}) follows.
\end{proof}

Lemma \ref{strict} has an important role in what follows. It was proved in the real case for smooth strictly convex functions in Lemma 2.2 of \cite{Slod}. In Lemma 10 of \cite{PS}, the authors translate it into the setting of plurisubharmonic functions on complex manifolds. Applying the same proof with necessary modifications (by using the standing assumption (*) and Lemma \ref{strlem} in the relevant parts) also carries it over to the $m$-subharmonic functions.

\begin{lemma} \label{strict}
Suppose we are given a compact subset $L$ of a K\"{a}hler manifold $M$ and a bounded smooth strictly $m$-subharmonic function $\rho$ which is defined in a neighborhood $U$ of $L$. Let $v$ be an $m$-subharmonic function on $U$. Suppose further that there is a non-empty compact set $K \subset L$ such that \begin{equation}\label{mas}\max_{L}{v}> \max_{K}{v}.\end{equation}Then there exist a point $z_{0} \in L \backslash K$, a neighborhood $V$ of $z_{0}$ and a smooth strictly $m$-subharmonic function $u$ on $V$ such that $u(z_{0})=0$ ; and  whenever $z_{0} \neq z \in L \cap V$, we have $u(z)<0$.
\end{lemma}

\begin{remark}\label{mlocc}
Following the ideas of the proof of Proposition 2.3 in \cite{Slod}, we can show, by using our standing assumption(*) and  Lemma \ref{strict} where necessary, the equivalence of (ii), (iv) and (v) in Proposition 2.3 for $m$-subharmonic functions. As a result of this, we see that $E$ is an $m$-local maximum set if and only if for any relatively compact open set $V\subset M$, $E\cap V$ is $m$-local maximum in $V$. In fact, as was obtained in Proposition 2.3 of \cite{Slod}, a similar version (maybe not necessarily using all the five items there) of local maximum sets for subharmonic functions can also be obtained due to Lemma 2.2 in \cite{Slod}. For other interesting details as to the local maximum property and its different types, we refer the reader to \cite{Slod}, \cite{HST}, \cite{HST2} and references therein.
\end{remark}

Next lemma is an easy adaptation of Proposition 1.8 in \cite{Slod2}. It can be proved exactly in the same way by using Lemma \ref{strict} and Proposition \ref{maxx} above.

\begin{lemma}\label{18}
Assume that we are given a continuous $m$-maximal function $u$ on $M$ and $v$ a $m$-subharmonic function on $M$ with $u(w_{0})=v(w_{0})$ and $v(w) \leq u(w)$ for all $w\in M$. Then $$G=\{z\in U: u(z)=v(z)\}$$ is an $m$-local maximum set.
\end{lemma}

By using Lemma \ref{strict} and arguments in the proof of Lemma 3.3 in \cite{MST}, one can prove the following lemma  that says that a level set of an $m$-subharmonic function inside an $m$-local maximum set is also $m$-local maximum.
\begin{lemma}\label{inmloc}
Let $X$ be a local $m$-maximum set in a K\"{a}hler manifold $M$, and $\varphi$
an $m$-subharmonic function defined in a neighborhood $U$ of $X$. Assume that $\varphi|_{X}$
reaches its absolute maximum value at some point $w_{0}$. Then the set $G= \{z \in X : \varphi(z) = \varphi(w_{0})\}$ has the $m$-local maximum property.

\end{lemma}

\begin{definition}\label{minim}
We call a function \textit{m-minimal function for the $m$-core $\mathbf{c}_{m}(M)$} if it is bounded from above continuous $m$-subharmonic function on $M$ which is also strictly $m$-subharmonic on $M \backslash \mathbf{c}_{m}(M)$.
\end{definition}

It is not difficult to see that, under the condition $\mathbf{c}_{m}(M) \neq M$, $m$-minimal functions do exist. Indeed, let us take an element $z_{q}\in M \backslash \mathbf{c}_{m}(M)$ and an open neighborhood $U_{q} \subset M \backslash \mathbf{c}_{m}(M)$ of $z_{q}$. Then there is $\psi_{q}\in \mathrm{SH}^{cb}_{m}(M)$ which is strictly $m$-subharmonic on $U_{q}$. Since $M$ is countable at infinity, that is, $M \backslash \mathbf{c}_{m}(M) = \bigcup_{q=1}^{\infty}{U_{q}}$, by the existence of  $\psi_{q} \in \mathrm{SH}^{cb}_{m}(M)$ being strictly $m$-subharmonic on $U_{q}$, one can define \begin{equation}\label{stpsh}\psi=\sum_{q=1}^{\infty}{\epsilon_{q} \psi_{q}}\end{equation} for a suitably chosen positive numbers $\{\epsilon_{q}\}$ such that (\ref{stpsh}) is uniformly convergent on compact subsets of $M$ and is an upper bounded $m$-subharmonic function on $M$. Then it follows from the compact-open topology of the Frechet space $\mathcal{C}(M)$ that $\psi$ is continuous on $M$. Lemma \ref{strlem} gives that it is strictly $m$-subharmonic on $\bigcup_{q=1}^{\infty}{U_{q}} = M \backslash \mathbf{c}_{m}(M)$ also.

It is important to note here that if $\mathbf{c}_{m}(M)= M$, then according to the $m$-subharmonic analogue of Corollary 5 (which can be done by using the same arguments there without any difficulty) from \cite{PS}, every $u\in \mathrm{SH}^{cb}_{m}(M)$ becomes $m$-maximal.

By mimicking the proof of Lemma 4.5 in \cite{Slod2} and using Lemma \ref{18}, Remark \ref{mlocc}, the standing local approximation condition(*) and Lemma \ref{maxloc} where necessary, we have the following theorem which gives us the $m$-pseudoconcavity of $\mathbf{c}_{m}(M)$.

\begin{theorem} \label{corepseudo}
Let $\varphi$ be a $m$-minimal function for $\mathbf{c}_{m}(M)$ and $B$ be a ball in $M$ that intersects with $\mathbf{c}_{m}(M)$. If $v: \overbar{B} \rightarrow \mathbb{R}$ be a continuous function, maximal on $B$ and $v|_{\partial B}= \varphi|_{\partial B}$, then $$B \cap \mathbf{c}_{m}(M) = \{z\in B: \varphi(z)=v(z) \}.$$ Therefore, $\mathbf{c}_{m}(M)$ is $m$-pseudoconcave.
\end{theorem}

In the rest of this section, we shall prove our main theorem, that is Theorem \ref{decompi}, by an argument  used in the proof of Theorem 3.3 of \cite{Slod2} as we mentioned before. Proposition \ref{capp} is used without proof in \cite{Slod2}. We shall provide its proof since it is crucial for the proof of main theorem. To this end, we will require another elementary topological lemma pertaining to upper semicontinuous functions, which might be somewhere in the literature, however we couldn't locate any reference to it, so we will supply a proof for it as well.

\begin{lemma}\label{usctop}
Let $\mathcal{F}=\{ F_{\alpha}\}_{\alpha \in \Lambda}$ be a family of compact sets in a complex manifold $M$ that is closed under taking finite intersections.  Let  $h$  be an upper semicontinuous function on the union $\bigcup_{\alpha \in \Lambda}{F_{\alpha}}$. Let $b\in \mathbb{R}$  such that for every set $F_{\alpha}\in \mathcal{F}$, $\max_{F_{\alpha}}{h} \geq b$. Let $\Gamma= \bigcap_{\alpha \in \Lambda}{F_{\alpha}}$. Then \begin{equation}\label{gamma}\max_{\Gamma}{h} \geq b\end{equation}in case $\Gamma$ is non-empty.
\end{lemma}

\begin{proof}
We make some observations, first of all, $\mathcal{F}=\{F_{I}=\cap_{\beta \in I}{F_{\beta}}: I \subset \Lambda, \,\text{where $I$ is any finite subset}\}$. From this, we also have $$\bigcap_{I \subset \Lambda}{F_{I}}=\bigcap_{\alpha \in \Lambda}{F_{\alpha}}=\Gamma.$$Since each $F_{I}\in \mathcal{F}$ by the property of the family $\mathcal{F}$, one gets $$\max_{F_{I}}{h} \geq b.$$

Let us return to the proof of (\ref{gamma}). Assume the contrary that $$\max_{\Gamma}{h}<b.$$Then the following set $$V=\{z\in \bigcup_{\alpha \in \Lambda}{F_{\alpha}}: h(z)< b\} $$ is an open set in $\bigcup_{\alpha \in \Lambda}{F_{\alpha}}$ (in the subspace topology induced from $M$) by the upper semicontinuity of $h$. This also gives that $\Gamma \subset V$.

Let $\alpha_{0}\in \Lambda$. Then one has  $$F_{\alpha_{0}} \subset (\bigcup_{\alpha_{0}\neq \alpha}{F^{c}_{\alpha}}) \cup \Gamma \subset (\bigcup_{\alpha_{0} \neq \alpha}{F_{\alpha}}^{c}) \cup V. $$ Since $F_{\alpha_{0}}$ is compact, there is a finite set $J \subset \Lambda$ such that $$F_{\alpha_{0}} \subset (\bigcup_{\alpha \in J}{F^{c}_{\alpha}}) \cup V,$$ which gives us $$\bigcap_{\alpha \in J \cup \{\alpha_{0}\}}{F_{\alpha}} \subset V.$$ Therefore we have, by writing  $J'=J \cup \{\alpha_{0}\}$, that $\Gamma \subset \cap_{\alpha \in J'}{F_{\alpha}} \subset V$. This concludes, by the property of the family $\mathcal{F}$, that $$ b \leq \max_{\bigcap_{\alpha\in J'}{F_{\alpha}}}{h} \leq \max_{V}{h} < b,$$which is a contradiction.

\end{proof}

\begin{proposition}\label{capp}
Let $G_{\alpha}, \,\alpha \in \Lambda$,\, be closed subsets of $M$. Suppose that for every finite subset $\{\alpha_{1}, \ldots, \alpha_{n}\} \subset \Lambda$, the intersection $\cap_{j=1}{G_{\alpha_{j}}}$ is an $m$-local maximum set. Then the set $$G=\bigcap_{\alpha \in \Lambda}{G_{\alpha}}$$ is also an $m$-local maximum set whenever it is non-empty.
\end{proposition}

\begin{proof}
Fix any $w \in G$ and a relatively compact open neighborhood $V$ of $w$. Let $h$ be an $m$-subharmonic function in $V$ with $h(w)=b$. Write $K= \partial V \cap G$ . We need to show that $\max_{K}{h}\geq b$.
Let now $L$ be the family of all finite intersections of sets in $\{G_{\alpha}\}$ with $\partial V$, such as $G_{t_{1}} \cap \ldots \cap G_{t_{n}} \cap \partial V$. Since finite intersections are local $m$-maximum sets by our hypothesis, we have $\max_{L_{\beta}}{h} \geq b$ for all $L_{\beta} \in L$. Also $L$ is preserved by finite intersections. Since $G\cap \partial V$ is the intersection of all sets $L_{\nu}$ in $L$, by Lemma \ref{usctop}, we obtain that $\max_{K}{h}\geq b$, which is what we needed.

\end{proof}

Next step is to get the $m$-subharmonic version of Lemma 3.5 in \cite{Slod2}. We cannot use the first part of the proof of Lemma 3.5 directly in our setting because the construction there was considered for bounded (from above and below) functions. The key point here is to utilize Lemma 1.14 in \cite{Slod2}. This lemma gives us the following: For every $b\in \mathbb{R}^{n+1}$ and for every $\epsilon>0$, one can construct a smooth convex function on $\mathbb{R}^{n+1}$ such that \textit{a)} for every $x\in \mathbb{R}^{n+1} \backslash \{b\}$, $v(x_{0}, \ldots, x_{n})> x_{0}+ \ldots +x_{n}$ ; \textit{b)} $v(b)=b_{0}+ \ldots + b_{n}$; \textit{c)} $\frac{\partial v}{\partial x_{j}}>0$ for every $j=0, \ldots, n$; \textit{d)} $v$ is of linear growth, i.e., $|v(x)|\leq c_{0}+ (1+ \epsilon)|x|$ for any $x\in \mathbb{R}^{n}$. Let $\varphi_{0}, \varphi_{1}, \ldots, \varphi_{n}$ be bounded from above continuous $m$-subharmonic on $M$ with $\varphi_{0}$ strictly $m$-subharmonic on $M \backslash \mathbf{c}_{m}(M)$, i.e., it is an $m$-minimal function. As in the proof of Lemma 4.8 in \cite{Slod2}, by using the function $v$ given above, we can define $\nu: M \rightarrow \mathbb{R}$ by $\nu(z)= v(\varphi_{0}, \ldots, \varphi_{n})$. This function is $m$-minimal function for the $m$-core $\mathbf{c}_{m}(M)$, that is to say, it is upper bounded, continuous $m$-subharmonic function and it is also strictly $m$-subharmonic on $M \backslash \mathbf{c}_{m}(M)$ because $\frac{\partial v}{\partial x_{j}}>0$  and $\varphi_{0}$ is strictly $m$-subharmonic on $M \backslash \mathbf{c}_{m}(M)$. Now define $\varphi= \varphi_{0}+ \ldots + \varphi_{n}$. Since, by $(b)$ above, $v(b)=b_{1}+ \ldots + b_{n}$, we have $\nu(z) \geq \varphi(z)$ for any $z\in M$ and therefore $$\{z\in M: \nu(z)=\varphi(z)\} = \{z\in M: \varphi_{j}(z)=b_{j},\,\,j=0, 1, \ldots, n \}.$$ By using the maximality argument as in the proof of Lemma 3.5 in \cite{Slod2} and the applications of  Theorem \ref{corepseudo} and Lemma \ref{18}, we have

\begin{lemma}\label{finitary}
Given that the functions $\varphi_{0}, \varphi_{1}, \ldots, \varphi_{n}$ are bounded from above continuous $m$-subharmonic on $M$ with $\varphi_{0}$ strictly $m$-subharmonic on $M \backslash \mathbf{c}_{m}(M)$, i.e., it is an $m$-minimal function. Let $b_{0}, \ldots, b_{n} \in \mathbb{R}$. Consider the set $$G=\mathbf{c}_{m}(M) \cap \bigcap_{l=0}^{n}{\{z\in M: \varphi_{l}(z)=b_{l}\}}.$$If $G$ is non-empty, then it is an $m$-local maximum set.
\end{lemma}

We define a relation now as follows: $z\sim w$ if $\varphi(z)=\varphi(w)$ for all $\varphi \in \mathrm{SH}^{cb}_{m}(M)$. It is easy to verify that this relation is an equivalence relation.

Note that, according to this equivalence relation, the equivalence class $[w]$ for any $w\in M$ is \begin{equation}\label{den} (A^{cb}_{m}(w))_{e}=\bigcap_{s\in S}{\{z\in M: \varphi_{s}(z)=\varphi_{s}(w)\}}.\end{equation}
Let us show that an equivalence class of an element which is not in the $m$-core reduces to a singleton.
\begin{lemma}\label{eqq}
If $w\notin \mathbf{c}_{m}(M)$, then $(A^{cb}_{m}(w))_{e}=\{w\}$.
\end{lemma}

\begin{proof}
Let us assume for contradiction that there is another element $z\in (A^{cb}_{m}(w))_{e}$ with $z\neq w$. By hypothesis, there is $\phi \in \mathrm{SH}^{cb}_{m}(M)$ being strictly $m$-subharmonic around an open neighborhood $V$ of $w$. By the definition of strict $m$-subharmonicity, for a function $\alpha \in \mathcal{C}^{\infty}_{0}(M)$ whose compact support $\supp{(\alpha)}$ is contained in $V$ with $\alpha(z)=1$ and $\alpha(w)=0$ (in other words, $w \notin \supp{(\alpha)}$), we can find $\epsilon_{0}>0$ such that $\phi+ \epsilon \alpha$ is $m$-subharmonic in $V$ for every $\epsilon \in (-\epsilon_{0}, \epsilon_{0})$. Fix $\epsilon>0$. Now for any $z' \notin \supp{(\alpha)}$, we get $(\phi+ \epsilon \alpha)(z')=\phi(z')$, namely $\phi+ \epsilon \alpha$ is $m$-subharmonic in $M \backslash \supp{(\alpha)}$. Therefore, $\phi+ \epsilon \alpha \in \mathrm{SH}^{cb}_{m}(M)$. Now $\phi(z)+ \epsilon=(\phi+ \epsilon \alpha)(z)=(\phi+ \epsilon \alpha)(w)=\phi(w)$, contradicting  $\phi(z)=\phi(w)$.
\end{proof}

Next corollary tells us that any equivalence class of an element in the $m$-core will be completely included in the $m$-core, that is,

\begin{corollary}\label{ewq}
For any $w\in \mathbf{c}_{m}(M)$, $(A^{cb}_{m}(w))_{e} \subset \mathbf{c}_{m}(M)$.
\end{corollary}
\begin{proof}
Suppose otherwise that there is an element $z\in (A^{cb}_{m}(w))_{e} \backslash \mathbf{c}_{m}(M)$. It is obvious that $z\neq w$ since $w\in \mathbf{c}_{m}(M)$. By the symmetry of the sets $(A^{cb}_{m}(w))_{e}$ and $(A^{cb}_{m}(z))_{e}$, we have $w\in (A^{cb}_{m}(z))_{e}$, which means that $(A^{cb}_{m}(z))_{e} \neq \{z\}$. Lemma \ref{eqq} yields that $z\in \mathbf{c}_{m}(M)$, contradiction.
\end{proof}
\begin{proof}[Proof of Theorem \ref{decompi}] Consider the family $\{\varphi_{s}: s\in S\}$ of all bounded from above continuous
$m$-subharmonic functions on $M$. Take $w\in \mathbf{c}_{m}(M)$. Let $G$ be the component of $\mathbf{c}_{m}(M)$ containing $w$. Write $\varphi_{s}(w)=b_{s}$. Then we have
$$G=\mathbf{c}_{m}(M) \cap (\bigcap_{s\in S}{\{z\in M: \varphi_{s}(z)=b_{s}\}})=\mathbf{c}_{m}(M) \cap (A^{cb}_{m}(w))_{e} $$by the relation (\ref{den}). Now by Corollary \ref{ewq}, one has $$G=\bigcap_{s\in S}{\{z\in M: \varphi_{s}(z)=b_{s}\}}=(A^{cb}_{m}(w))_{e}.$$

Since the sets $(A^{cb}_{m}(w))_{e}$ are equivalence classes, they are either disjoint or identical, so we have the following disjoint union of equivalence classes of points of the $m$-core $$\mathbf{c}_{m}(M)= \bigsqcup_{w\in \mathbf{c}_{m}(M)}{(A^{cb}_{m}(w))_{e}}.$$
What remains to see is that each $(A^{cb}_{m}(w))_{e}$ is an $m$-pseudoconcave set. For this, let us take $s_{0}\in S$ so that $\varphi_{s_{0}}$ is an $m$-minimal function. Take any finite set $I \subset S$ with $s_{0}\in I$ and define $$G_{I}= \mathbf{c}_{m}(M)\cap (\bigcap_{j\in I}{\{z\in M: \varphi_{s_{j}}(z)=b_{s_{j}}\}}).$$ Lemma \ref{finitary} implies that all sets of the form $G_{I}$ are $m$-local maximum sets. Since $G= \bigcap_{I\subset S}{G_{I}}$, by Proposition \ref{capp} and Lemma \ref{maxloc}, we obtain that $G$ is $m$-pseudoconcave.
\end{proof}
Combining one of the main ingredients in the proof of Theorem \ref{decompi} saying that $(A^{cb}_{m}(w))_{e}$ is $m$-pseudoconcave for $w\in \mathbf{c}_{m}(M)$ with Lemma \ref{eqq}, the corollary below is immediate

\begin{corollary}
$w\in \mathbf{c}_{m}(M)$ if and only if\, $(A^{cb}_{m}(w))_{e}\neq \{w\}$.
\end{corollary}

We now give some examples of domains with cores of various kinds.

\begin{example} Let $M$ be a Stein manifold. Since $M$ can be holomorphically embedded into $\mathbb C^N$ for some $N$, any $m-$subharmonic  function can be made strictly $m-$subharmonic by adding $\epsilon |z|^2$ to it, where $z=(z_1,\dots ,z_N)$ are the local coordinates. Thus $\mathbf{c}_{m}(M)=\emptyset.$
\end{example}

\begin{example}\label{ex2} Let $M=\{(z,w)\in \mathbb C^2: \log|z-w|+|z|^2+|w|^2<0\}$ and $L=\{(z,z): z\in \mathbb C \}\subset M.$ Any upper bounded $1-$subharmonic function is constant on $L$. Thus $\mathbf{c}_{1}(M)\supset L$. $\psi(z,w)=\log |z-w|+|z|^2+|w|^2$ is continuous, bounded above, strictly $1$-subharmonic on $M\setminus L$. Thus $\mathbf{c}_{1}(M)\subset L$ and hence $\mathbf{c}_{1}(M)=L$. Observe that $L$ is connected.
\end{example}

\begin{example} Let $M=\{(z,w)\in \mathbb C^2:\log|z|+\log|z-1|+|z|^2+|w|^2<0  \}$.  Then $$\mathbf{c}_{1}(M)=\{z=0\}\cup \{z=1\}$$ is not connected.

\end{example}

\begin{example}Let $M=\mathbb D\times \mathbb C.$ Any upper bounded continuous $1$-subharmonic function is constant on the line $\{z=z_0\}\subset M$ for any $z_0\in \mathbb D$. Thus $\mathbf{c}_{1}(M)=M$.

\end{example}

\section{a characterization of pluriharmonic functions}\label{four}

In this last section, we will be interested in so-called $m$-harmonic functions in the context of $m$-subharmonic functions. We understand that, on a non-compact K\"{a}hler manifold, pluriharmonic functions and $m$-harmonic functions($m\geq 2$) are the same. We start with a definition. Let $M$ be a (non-compact) K\"{a}hler manifold with a K\"{a}hler form $\omega$ on it and $\Omega \subset M$ a domain. In this section we study pluriharmonic functions in terms of so-called $m$-harmonic functions.  We call $v\in \mathcal{C}^{2}(\Omega)$ \, \textsl{m-harmonic} in case \begin{equation}\label{mhar}
(dd^{c}v)^k \wedge \omega^{n-k}= 0 \end{equation} holds for $k=1, 2, \ldots, m$. By definition, any $m-$harmonic function is $(m-k)-$harmonic for any $k=1,\dots , m-1$. Recall that $v$ is pluriharmonic if $dd^c v=0$. Thus pluriharmonic functions are $m-$harmonic for any $m=1,\dots ,n.$

\begin{remark}
A related concept was defined with the name $k$-harmonic in \cite{LN} for mappings from balls of dimension $m$ into balls of dimension $n$. Our definition of $m$-harmonicity is different than this definition. The $m$-harmonicity that we introduced in this section satisfies the geometric condition given in \cite{LN} not only on $k$-dimensional subspace but also on all complex subspaces of dimension lower than $k$. Another important point here is that those functions which verify this geometric property on all possible lower dimensional subspaces form a strictly bigger class than $m$-harmonic functions that we have defined above. This means that we cannot characterize $m$-harmonic functions and more generally $m$-subharmonic functions by this geometric feature.  For more details, see \cite{Din} and \cite{Sa}.
\end{remark}

In order to prove the equivalence of $m$-harmonic and pluriharmonic functions generally on a K\"{a}hler manifold, first we pay attention to the local case, namely, $M=\mathbb{C}^{n}$, where $n\geq 2$, and show that this equivalence is true in $\mathbb{C}^n$. We now record some information that we need about complex $k$-Hessian operators. For more details and background, we refer the reader to the papers mentioned in the section 1, for example, \cite{Bl} or \cite{SA} are comprehensive.

Given that $v\in \mathcal{C}^{2}(\Omega)$,  $$dd^{c}v=\frac{i}{2}\sum_{k,l}{\frac{\partial^{2} v}{\partial z_{k} \partial \overbar{z_{l}}}dz_{k} \wedge d\overbar{z_{l}}}$$ is an Hermitian quadratic form. Using a suitable unitary transformation, we can turn the above second order differential expression into a diagonal form as follows: \begin{equation}\label{diagon}dd^{c}v= \frac{i}{2}(\sum_{l=1}^{n}\nu_{l} dz_{l} \wedge d\overbar{z_{l}}),\end{equation} where $\nu_{1}, \ldots, \nu_{n}$ are the (real) eigenvalues of the Hermitian matrix $[\partial ^{2} v/ \partial z_{k} \partial \overbar{z_{l}}]^{n}_{k,l=1}$.

 For $m=1$, we get the harmonic functions on $\mathbb{C}^{n}\simeq \mathbb{R}^{2n}$. Clearly any $m$-harmonic function is $m$-subharmonic.

In $\mathbb{C}^{n}, n \geq 2$, assuming as above that $v$ is at least of $\mathcal{C}^{2}$, we have the following explicit form \begin{equation}\label{hess}
(dd^{c}v)^k \wedge \omega^{n-k}=4^{n} k! (n-k)! H_{k}(v)d\lambda, \,\,\,k=1, 2, \ldots, m,
\end{equation}where $\omega=dd^{c}|z|^{2}$ is the fundamental K\"{a}hler form on $\mathbb{C}^{n}$, $d\lambda$ is the volume form on $\mathbb{C}^n$ and \begin{equation}\label{symm}
H_{k}(v)=\sum_{1\leq j_{1}< \ldots < j_{k} \leq n}{\nu_{j_{1}} \ldots \nu_{j_{k}}}
\end{equation} is the elementary complex $k$-Hessian operator. We know any pluriharmonic function is also $m$-harmonic as remarked above, but the converse might not be true in general, because for $m=1$, that is to say, for harmonic functions on $\mathbb{C}^{n} \simeq \mathbb{R}^{2n}$, we have $$dd^{c}v \wedge \omega^{n-1}=0$$ implies, by (\ref{hess}) and (\ref{symm}), that $H_{1}(v)= \nu_{1}+ \ldots + \nu_{n}=0,$ but from this equation, not all the eigenvalues $\nu_{1}, \ldots, \nu_{n}$ have to be zero and it is well-known that the class of harmonic functions properly contains that of pluriharmonic functions. However, if we consider $m\geq 2$, then the converse becomes true. Let us see this now. Since the elementary complex $k$-Hessian operator $$H_{k}(v)=\sum_{1\leq j_{1}< \ldots < j_{k} \leq n}{\nu_{j_{1}} \ldots \nu_{j_{k}}}, \,\,\,k=1, 2, \ldots, m$$ is an elementary symmetric polynomial in the variables $\nu_{1}, \ldots, \nu_{n}$, it will suffice to consider the cases $m=1$ and $m=2$. For these two cases, we get $$H_{1}(v)=\sum_{k=1}^{n}{\nu_{k}}=0,\,\,\,H_{2}(v)=\sum_{1\leq k<l\leq n}{\nu_{k}\nu_{l}}=0.$$According to this last information, $\nu_{1}= \ldots =\nu_{n}=0$, which yields, by using the relation (\ref{diagon}), $dd^{c}v=0$, i.e. $v$ is pluriharmonic. What seems to be interesting in the above discussion is that we do not have to look at all of $m$ equations that come from the $m$-harmonicity relation (\ref{mhar}), the first two cases give the pluriharmonicity. Hence we have proved

\begin{theorem}\label{mheq}Let $\Omega \subset \mathbb{C}^{n}, \,\,n\geq 2,$ be a domain and $v\in \mathcal C^2(\Omega).$ Then $v$ is $2$-harmonic if and only if $v$ is pluriharmonic.
\end{theorem}

Working with local coordinates, we can generalize the above result for K\"ahler manifolds.

\begin{corollary}
 Theorem \ref{mheq} holds for K\"ahler manifolds of dimension at least $2$.
\end{corollary}

\begin{proof} Let $(M,\tilde{\omega})$ be a K\"ahler manifold of dimension at least $2$ with local coordinates $z=(z_1,\dots, z_n)$ near $p\in M$. For simplicity, we take $p=0$ and $\omega=dd^c|z|^2$ denotes the standard K\"ahler form. Near $0$, $\tilde{\omega}$ osculates to order two to standard form $\omega$. That is, $\tilde \omega=\omega +O(|z|^2)\alpha$, where $\alpha$ is the $(1,1)$ form $\sum_{i,j}dz_i\wedge d\bar z_{j}.$ Then we have $$\tilde \omega^{n-m}=\omega^{n-m}+\sum_{j=0}^{n-m-1}{{n-m}\choose j}\omega^j\wedge (O(|z|^2)\alpha)^{n-m-j}=\omega^{n-m}+O(|z|^2)\beta, $$ where $\beta$ is an $(n-m,n-m)$ form.  Hence we get $$(dd^c v)^m\wedge \tilde \omega^{n-m}=4^nm!(n-m)!H_m(v)d\lambda +O(|z|^2)(dd^c v)^m\wedge \beta $$ where $H_m(v)$ is defined in (\ref{symm}). If $H_m(v)>0$ at $0$, then there exists a neighborhood $V$ of $0$ such that   $4^nm!(n-m)!H_m(v)>\delta$ for some $\delta>0$. Take a smaller neighborhood $V'\subset V$ of $0$ such that $O(|z|^2)(dd^c v)^m\wedge \beta\leq \frac{ \delta}{2}.$ Then on $V'$, $ (dd^c v)^m\wedge \tilde \omega^{n-m}>\frac{\delta}{2}$. Accordingly, if $(dd^c v)^m\wedge \tilde \omega^{n-m}=0$ then  $H_m(v)=0$.  It follows from the same argument in the proof of Theorem \ref{mheq} that any $m$-harmonic function $(m\geq 2)$ is pluriharmonic on $M$.
\end{proof}

\textbf{Acknowledgments.} We are thankful to Zbigniew Slodkowski for providing several materials from his paper (Pseudoconcave decompositions in complex manifolds) and for answering our questions and useful discussions.

 The first and second author of this paper are supported by the T\"UB\.ITAK  2518 Proj. No. 119N642. The third author is supported by T\"UB\.ITAK 2518 Proj. No. 119N642 and  T\"UB\.ITAK 3501 Proj. No. 120F084.

We gratefully thank the anonymous referee for improving the exposition of the paper.

\end{document}